\documentclass{article}
\usepackage{amsmath}
\usepackage{amsfonts}   
\usepackage{amssymb}

\newcommand{\mlabel}[1]{\label{#1}
}

\newcommand{\seq}{\begin{equation}}                 
\newcommand{\eeq}[1]{\label{#1}\end{equation}
    }

\newcommand{\pf}{ \par \vspace{1ex} \noindent {\sc Proof} \hspace{2mm}}
\newcommand{\epf}{$ \quad \Box$ \par \vspace{1ex}}

\newtheorem{Theorem}{Theorem}[section]
\newcommand{\sthm}{\begin{Theorem}}         
\newcommand{\ethm}{\end{Theorem}}           

\newtheorem{Corollary}[Theorem]{Corollary}
\newcommand{\scor}{\begin{Corollary}}       
\newcommand{\ecor}{\end{Corollary}}         

\newtheorem{Lemma}[Theorem]{Lemma}
\newcommand{\slm}{\begin{Lemma}}            
\newcommand{\elm}{\end{Lemma}}              

\newtheorem{Remark}[Theorem]{Remark}
\newcommand{\srmark}{\begin{Remark}\rm}        
\newcommand{\ermark}{\end{Remark}}             
\newtheorem{Example}[Theorem]{\sc Example}
\newcommand{\sex}{\begin{Example}\rm}        
\newcommand{\eex}{\end{Example}}             

\newcommand{\seql}{\begin{eqnarray*}}       
\newcommand{\eeql}{\end{eqnarray*}}

\newcommand{\smlist}[1]{\begin{list}           
                      {(#1{zzcount})}{\usecounter{zzcount}}}
\newcommand{\elist}{\end{list}}

\newcommand{\for}{\quad \mbox{\textrm{for}} \quad}

\newcommand{\e}{\varepsilon}

\newcommand{\G}{\Gamma}
\renewcommand{\l}{\lambda}

\newcommand{\var}{\varphi}
\newcommand{\s}{\sigma}

\renewcommand{\t}{\tau}
\renewcommand{\O}{\Omega}

\newcommand{\vect}[1]{\mathbf{#1}}
\providecommand{\abs}[1]{\lvert#1\rvert}
\providecommand{\norm}[1]{\lVert#1\rVert}

\begin{document}
\pagenumbering{arabic}
\title{\bf{\small{\uppercase{Convex Solutions of systems of Monge-Amp\`ere equations}}}\footnotetext{This manuscript is available at arxiv:1007.3013}}

\author{Haiyan Wang\\\\Division of Mathematical and Natural Sciences \\ Arizona State University\\Phoenix, AZ 85069-7100, U.S.A.\\E-mail: wangh@asu.edu }

\date{}

\maketitle

\begin{abstract}
The existence and multiplicity and nonexistence of nontrivial radial convex solutions  of systems of Monge-Amp\`ere equations
are established with superlinearity or sublinearity assumptions for an appropriately chosen parameter. The proof of the results
is based on a fixed point theorem in a cone.
\end{abstract}


\section{Introduction and main results}
In this paper we consider the existence,  nonexistence
of convex solutions for the boundary values problem

\begin{equation}\label{eq3}
\left\{ \begin{array}{llll}
\Big(\big(u_1'(r)\big)^N\Big)' &= \lambda N r^{N-1} f^1(-u_1,...,-u_n),\;\; 0<r<1 \\
... \\
\Big(\big(u_n'(r)\big)^N\Big)' &= \lambda N r^{N-1}  f^n(-u_1,...,-u_n),\;\; 0 < r <1,\\
 u_i'(0)=u_i(1) &= 0 \;\; , i=1,...,n,
\end{array} \right.
\end{equation}
where $N \geq 1$.  Such a problem arises in the study of the existence of
convex radial solutions to the Dirichlet problem for the system of the Monge-Amp\`ere equations
\begin{equation}\label{eq1}
\left\{ \begin{array}{llll}
\text{det}(D^2u_1) &= \lambda f^1(-u_1,...,-u_n)\;\; \text{in} \;\; B  \\
... \\
\text{det} (D^2u_n) &= \lambda   f^n(-u_1,...,-u_n)\;\; \text{in} \;\; B,\\
 u_i = 0 \;\; &\text{on} \;\; \partial B \;, i=1,...,n,
\end{array} \right.
\end{equation}
where $D^2 u_i = (\frac{\partial u_i}{\partial x_i \partial x_j})$ is the Hessian matrix of $u_i$, $B=\{x \in \mathbb{R}^N: |x| < 1 \}$.
Monge-Amp\`ere equations arise from Differential Geometry and optimization and mass-transfer problems.

It was  shown (e.g.
in \cite{Zhang2009}) that any convex solution of (\ref{eq-n=1}) must be radially symmetric for some special function $f$. Therefore it is reasonable to look
for convex radial solutions of (\ref{eq1}). For radial solution $u_i(r)$ with $r=\sqrt{\sum_1^N x_i^2}$, the Monge-Amp\`ere operator simply becomes
\begin{equation}\label{detOP}
\text{det} (D^2u_i)=\frac{(u_i')^{N-1}u_i''}{r^{N-1}}=\frac{1}{Nr^{N-1}}((u_i')^{N})'
\end{equation}
and then  (\ref{eq1}) can be easily transformed into (\ref{eq3}). (\ref{detOP}) is frequently used in the literature, see
the references below on radial solutions of the Monge-Amp\`ere equations and  others, e.g., Caffarelli and Li \cite{CAFFARELLILI2003}.
It can be derived from the fact that the Monge-Amp\`ere operator is rotationally invariant, see, for example, Goncalves and  Santos \cite[Appendix A.2]{Goncalves2005}.

Much attention has been focused on the study of the problem with a single equation, see e.g. \cite{CAFFARELLILI1984,CAFFARELLILI2003,Gutierrez2000}. When $n=1,$ (\ref{eq1}) reduces to
\begin{equation}\label{eq-n=1}
\begin{split}
\text{det} (D^2u)  &= \l f(-u) \;\; \text{in} \;\; B \\
u &=0 \;\; \text{on} \;\; \partial B.
\end{split}
\end{equation}
Lions \cite{LIONS1985} obtained existence results for (\ref{eq-n=1}) in general domains in $\mathbb{R}^N$
where the particular function $f(u)=u^N$ acts like a ``linear" term to the fully nonlinear
operation $det(D^2u)$. Kutev \cite{KUTEV1988} investigated the existence of strictly convex radial solutions of (\ref{eq-n=1}) with $f(-u)=(-u)^p$.
The author \cite{Wang2006}, and Hu and the author \cite{HU2006} showed that the existence, multiplicity and nonexistence of
convex radial solutions of (\ref{eq-n=1}) can be characterized by the asymptotic behaviors of the quotient $\frac{f(u)}{u^N}$ at zero and infinity.
The existence of convex radial solutions of some special Monge-Amp\`ere equations can also be found in \cite{DELANOE1985}.

Before stating our theorems, we make following assumptions. Let $\mathbb{R}=(-\infty, \infty)$, $\mathbb{R}_+=[0, \infty)$, $\mathbb{R}_+^n=\underbrace{\mathbb{R}_+ \times ... \times \mathbb{R}_+}_{n}$ and
$\norm{\vect{u}}=\sum_{i=1}^n \abs{u_i}$ for $\vect{u}=(u_1,...,u_n) \in \mathbb{R}^n_+.$

\begin{enumerate}
  \item[] \begin{itemize}
   \item[(H1)]$f^i: \mathbb{R}_+^n \to \mathbb{R}_+$ is continuous.
   \item[(H2)]$ f^i(u_1,...,u_n) > 0$ for $\vect{u}=(u_1,...,u_n) \in \mathbb{R}_+^n$ and $\norm{\vect{u}} > 0,$ $i=1,...,n.$
   \end{itemize}
\end{enumerate}

We shall use the following notation, for $\vect{u} \in \mathbb{R}_+^n$, $i =1,...,n$
$$ \var(t)=t^N, \; \var^{-1}(t)=t^{\frac{1}{N}}, t \geq 0, \; f_0^i =\lim_{\norm{\vect{u}} \to 0} \frac{f^i(\vect{u})}{\var(\norm{\vect{u}})},\quad  f_{\infty}^i =\lim_{\norm{\vect{u}} \to \infty} \frac{f^i(\vect{u})}{\var(\norm{\vect{u}})}$$
\begin{equation}\label{notation}
\vect{f}_0=\sum_{i=1}^n f_0^i, \quad\vect{f}_{\infty}=\sum_{i=1}^n f_{\infty}^i.
\end{equation}
The notation $\vect{f}_0$ and $\vect{f}_{\infty}$ were introduced in the author \cite{Wang2003} to
define superlinearity and sublinearity at $0$ and $\infty$ for general systems of ordinary differential equations involving $p$-Laplaican.  They are analogous
to $f_0=\lim_{u \to 0} \frac{f(u)}{u^N}$ or $f_{\infty}=\lim_{u \to \infty} \frac{f(u)}{u^N}$ for scalar equations (e.g. \cite{Wang2006,Wang1994}).
Our arguments here are closely related to those in \cite{Wang2003} for the existence of positive solutions of systems of equations
involving $p$-Laplaican.  With the special form of the Monge-Amp\`ere operator,
we are able to make sharper estimates on the operator. Apparently the estimates can be further improved. The intervals of parameter $\lambda$ for
ensuring the existence
of convex solutions of (\ref{eq3}) are not necessarily optimal. We will address them
in the future.  Lemmas \ref{lm2}, \ref{lm6} are the same as in \cite{Wang2003}, and  proved here only for completeness.

Our main results are Theorems \ref{th1} and \ref{th2}. By a nontrivial convex solution $\vect{u}$ to (\ref{eq3}), we understand a nonzero vector-valued function $\vect{u}(r)=(u_1(r),...,u_n(r))$ (at least one component is not zero)
such that
$ u_i \in C^2[0,1]$ and  convex on $[0,1]$, $i=1,...,n$ and satisfies (\ref{eq3}). A nontrivial convex solution of (\ref{eq3}) is negative on [0,1).

\sthm\mlabel{th1} Assume \rm{(H1)} holds.\\
(a). If $\vect{f}_0 =0$ and $\vect{f}_{\infty}=\infty$, then for all $\l > 0$ (\ref{eq3}) (and (\ref{eq1}))  has
a nontrivial convex solution. \\
(b). If $\vect{f}_0 =\infty$ and $\vect{f}_{\infty}=0$, then for all $\l > 0$ (\ref{eq3}) (and (\ref{eq1})) has
a nontrivial convex solution.
\ethm
\sthm\mlabel{th2} Assume \rm{(H1)-(H2)} hold.\\
(a). If  $\vect{f}_0 =0$ or $\vect{f}_{\infty}=0$, then there exists a $\l_0 > 0$ such that
for all $ \l > \l_0$ (\ref{eq3}) (and (\ref{eq1}))  has a strictly convex solution. \\
(b). If $\vect{f}_0 =\infty$ or $\vect{f}_{\infty}=\infty$, then there exists a $\l_0 > 0$ such that
for all $ 0< \l < \l_0$ (\ref{eq3}) (and (\ref{eq1})) has a strictly convex solution. \\
(c). If $\vect{f}_0=\vect{f}_{\infty}=0$, then there exists a $\l_0 > 0$ such that
for all $\l > \l_0$ (\ref{eq3}) (and (\ref{eq1})) has two strictly convex solutions. \\
(d). If $\vect{f}_0 =\vect{f}_{\infty}=\infty$, then there exists a $\l_0 > 0$ such that
for all $0< \l < \l_0$ (\ref{eq3}) (and (\ref{eq1})) has two strictly convex solutions.\\
(e). If $\vect{f}_0  < \infty$ and $\vect{f}_{\infty} < \infty$, then there exists a $\l_0 > 0$ such that
for all $0< \l < \l_0$ (\ref{eq3}) (and (\ref{eq1})) has no strictly convex radial solution. \\
(f). If $\vect{f}_0 > 0$ and $\vect{f}_{\infty} > 0$, then there exists a $\l_0 > 0$ such that
for all $\l > \l_0$ (\ref{eq3}) (and (\ref{eq1})) has no strictly convex radial solution.\\
\ethm

\srmark\label{rem1}
When $n=1$, (\ref{eq1}) is reduced to the single equation (\ref{eq-n=1}). $\vect{f}_0$ and $\vect{f}_{\infty}$
become $f_0=\lim_{u \to 0} \frac{f(u)}{u^N}$ and $f_{\infty}=\lim_{u \to \infty} \frac{f(u)}{u^N}$ respectively. In this case, (H2) becomes
$f(u)>0$ for $u>0$. Thus Theorems \ref{th1} and
\ref{th2} cover  relevant results in \cite{Wang2006,HU2006}.   In \cite{Wang2009}, the author discusses
the existence of convex solutions of (\ref{eq3}) or (\ref{eq1}) when $n=2$ and $f^i$ only depends on one variable. It is also worthwhile to
note that the results in this paper do not cover those in \cite{Wang2009}.
\ermark

\srmark\label{rem2}
A nontrivial solution $\vect{u}(r)=(u_i(r))$ of (\ref{eq3}) has at least one nonzero component. Some of its components can be zero.
Then $\vect{v}(r)=(v_i(r))=(-u_i(r))$ is a nontrivial solution to (\ref{eq4}) or a positive fixed point of $\vect{T}_{\l}$ in (\ref{T_def}).
From the integral expression  (\ref{T_def}), we have for $r \in (0,1)$
$$v'_i (r)=-\big(\lambda \int_0^r N \t^{N-1}f^i\big(v_1(\t),...,v_n(\t))d\t\big)^{\frac{1}{N}}$$
Therefore, each component $u'_i(r)=-v'_i(r)$ is nondecreasing and $u_i=-v_i(r)$ is convex.

If we assume both (H1) and (H2) hold, it follows from
Lemma \ref{lm2} that for $s \in (0,1)$,   $\sum_{i=1}^n v_i(s)>0$ and
$$
\int^s_0 N  \t^{N-1}f^i\big(v_1(\t),...,v_n(\t)d\t >0.
$$
Thus (\ref{T_def}) implies that
\begin{equation*}
\begin{split}
&v''_i (r) \\
 &= -\frac{1}{N} \Big(\lambda\int^r_0 N  \t^{N-1}f^i\big(v_1(\t),...,v_n(\t)\big)d\t\Big)^{\frac{1}{N}-1}\big(\lambda N r^{N-1}f^i((v_1(r),...,v_n(r))\big)\\
& <0.
\end{split}
\end{equation*}
for $r \in (0,1)$. Thus if (H2) holds, each component $u_i$ must be strictly convex.
\ermark

We now give the following two examples to demonstrate a few  cases of the two theorems.  Examples for other cases  can be constructed  in the same way.


\textbf{Example 1 }
\begin{equation}\label{eq1-ex1}
\left\{ \begin{array}{lllll}
\text{det} (D^2u_1) &= \lambda (-u_1-u_2)^{p_1}\;\; \text{in} \;\; B  \\
\text{det} (D^2u_2) &= \lambda (-u_1-u_2)^{p_2}\;\; \text{in} \;\; B,\\
u_1=u_2 = 0 \;\; &\text{on} \;\; \partial B,
\end{array} \right.
\end{equation}
where $B=\{x \in \mathbb{R}^N: |x| < 1 \}$ as before. If $p_1,p_2>N$, then $\vect{f}_0 =0$ and $\vect{f}_{\infty}=\infty$.
Then (\ref{eq1-ex1}) has a strictly convex radial solution  for all $\l>0$ according to Theorem \ref{th1}. In the same way, If $0<p_1,p_2<N$,
then  $\vect{f}_0 =\infty$ and $\vect{f}_{\infty}=0$.  Then (\ref{eq1-ex1}) has a strictly convex radial solution  for all $\l>0$
according to Theorem \ref{th1}. On the other hand, if $0<p_1<N$ and $p_2>N$, then $\vect{f}_0 =\infty$ and $\vect{f}_{\infty}=\infty$.
According to Theorem \ref{th2}, (\ref{eq1-ex1}) has two strictly convex radial solutions for sufficiently small $\l>0$.

\textbf{Example 2 }
\begin{equation}\label{eq1-ex2}
\left\{ \begin{array}{lllll}
\text{det} (D^2u_1) &= \lambda e^{-u_1-u_2}\;\; \text{in} \;\; B  \\
\text{det} (D^2u_2) &= \lambda g(-u_1,-u_2)\;\; \text{in} \;\; B,\\
u_1=u_2 = 0 \;\; &\text{on} \;\; \partial B,
\end{array} \right.
\end{equation}
where $B=\{x \in \mathbb{R}^N: |x| < 1 \}$ as before $g$ is a any continuous function such that $\lim_{\norm{\vect{u}} \to 0}g$ and
$\lim_{\norm{\vect{u}} \to \infty}g$ (In fact, we may only need $\limsup,\liminf $ ) exist. Then $\vect{f}_0 =\infty$ and $\vect{f}_{\infty}=\infty$.
Then (\ref{eq1-ex2}) has no nontrivial convex radial solution for sufficiently large $\l>0$ and two nontrivial convex radial solutions for sufficiently small $\l>0$
according to Theorem \ref{th2}.


\section{Preliminaries}

We shall treat convex classical solutions of (\ref{eq3}), namely a vector-valued function $\vect{u}(t)$
of class $C^2[0, 1]$,  satisfying (\ref{eq3}). For the remaining sections, $t$ is often used as independent variable of functions and $r$ as radiuses of balls in the cone.  Now we treat positive concave classical solutions of (\ref{eq4}).
With a simple transformation $v_i=-u_i$ (\ref{eq3}) can be brought to the following equation

\begin{equation}\label{eq4}
\left\{ \begin{array}{llll}
\Big(\big(-v_1'(t)\big)^N\Big)' &= \lambda N t^{N-1} f^1(v_1,...,v_n),\;\; 0<t<1 \\
... \\
\Big(\big(-v_n'(t)\big)^N\Big)' &= \lambda N t^{N-1}  f^n(v_1,...,v_n),\;\; 0 < t <1,\\
 v_i'(0)=v_i(1) &= 0, \;\; i=1,...,n,
\end{array} \right.
\end{equation}
 We recall some concepts and conclusions of an operator in a cone. Let $X$
be a Banach space and $K$ be a closed, nonempty subset of $X$. $K$
is said to be a cone if $(i)$~$\alpha u+\beta v\in K $ for all
$u,v\in K$ and all $\alpha,\beta \geq 0$ and $(ii)$~$u,-u\in K$ imply
$u=0$. We shall use the following well-known fixed
point theorem to prove Theorems \ref{th1}, \ref{th2}.

\slm\mlabel{lm1} {\rm (\cite{DEIMING, GUOL, KRAS})} Let $X$ be a Banach
space and $K\ (\subset X)$ be a cone. Assume that $\Omega_1,\
\Omega_2$ are bounded open subsets of $X$ with $0 \in \Omega_1,\bar\Omega_1 \subset \Omega_2$, and let
$$
T: K \cap (\bar{\Omega}_2\setminus \Omega_1 ) \rightarrow K
$$
be completely continuous such that either
\begin{itemize}
\item[{\rm (i)}] $\| Tu \| \geq \| u \|,\ u\in K\cap \partial
     \Omega_1$ and $ \| Tu \| \leq \| u \|,\ u\in K\cap \partial
     \Omega_2$; or

\item[{\rm (ii)}] $\| Tu \| \leq \| u \|,\ u\in K\cap \partial
     \Omega_1$ and $\| Tu \| \geq \| u \|,\ u\in K\cap \partial
     \Omega_2$.
\end{itemize}
Then $T$ has a fixed point in $K \cap ( \bar \Omega_2 \backslash
     \Omega_1)$.

\elm

%
%

In order to apply Lemma \ref{lm1} to (\ref{eq1}), let $X$ be the Banach space \\$\underbrace{C[0,1] \times ... \times C[0,1]}_{n}$
and, for $\vect{v}=(v_1,...,v_n) \in X,$
$$\displaystyle{\norm{\vect{v}}= \sum_{i=1}^n \sup_{t\in[0,1]} \abs{v_i(t)}}.$$  For $\vect{v} \in X$ or $\mathbb{R}^n_+$, $\norm{\vect{v}}$ denotes the norm of $\vect{v}$ in $X$ or
$\mathbb{R}^n_+$, respectively.

Let $K$ be a cone in $X$ defined as
\begin{equation*}
\begin{split}
K  = \{&\vect{v}=(v_1,...,v_n) \in X: v_i(t)\geq 0,\; t \in [0,1], \; i=1,...,n,\\
&\rm{and} \; \min\limits_{\frac{1}{4}\leq t \leq \frac{3}{4}}\sum_{i=1}^n\text{$v_i(t)$} \geq \frac{1}{4} \norm{\vect{v}}\}.
\end{split}
\end{equation*}
For $r>0$ let
$$
\O_r  = \{\vect{v} \in K: \norm{\vect{v}} < r \}.
$$
Note that $\partial \O_r = \{\vect{v} \in K: \norm{\vect{v}}=r\}$.

Let $\vect{T}_{\l}: K \to X$ be a map with components $(T_{\l}^1,...,T_{\l}^n)$. We define $T_{\l}^i$, $i=1,...,n$, by
\begin{equation}\label{T_def}
T_{\l}^i\vect{v}(t) =\int^1_t \var^{-1}\Big(\l \int^s_0 N  \t^{N-1}f^i\big(v_1(\t),...,v_n(\t)\big)d\t\Big)ds, \;\; 0 \leq t \leq 1.
\end{equation}
Thus, if $\vect{v} \in K$ is a positive fixed point of $\vect{T}_{\l}$, then $-\vect{v}$ is a nontrivial
convex solution of (\ref{eq3}) or (\ref{eq1}). Conversely,
if $\vect{v}$ is a convex radial solution of (\ref{eq3}) or (\ref{eq1}), then $-\vect{v}$ is a fixed point of $\vect{T}_{\l}$ in $K$.

The following lemma is a standard result due to the concavity of $v$, see e.g. \cite{Wang2003}. We prove it here only for completeness.
\slm\mlabel{lm2}
Assume \rm{(H1)} hold. Let $v(t) \in C^1[0,1]$ for $ t \in [0,1]$. If $v(t) \geq 0 $ and $v'(t)$ is nonincreasing on $[0,1]$.
Then
$$
v(t) \geq \min\{t, 1-t\}||v||, \quad t \in [0, 1]
$$
where $||v||=\max_{t \in [0,1]}v(t).$ In particular, $$
\min_{\frac{1}{4} \leq t \leq \frac{3}{4}}v(t) \geq \frac{1}{4} ||v||.
$$
and if $v(0)=||v||$, then
$$
v(t) \geq (1-t)||v||, t \in [0,1].
$$

\elm
\pf Since $v'(t)$ is nonincreasing, we have for \mbox{$0 \leq t_0 < t < t_1 \leq 1,$}
$$
v(t)-v(t_0)=\int^t_{t_0}v'(s)ds \geq (t- t_0)v'(t)
$$
and
$$
v(t_1)-v(t)=\int^{t_1}_{t} v'(s)ds \leq (t_1 - t)v'(t),
$$
from which, we have
$$
v(t) \geq  \frac{(t_1-t)v(t_0) + (t-t_0)v(t_1)}{t_1-t_0}.
$$
Choosing $\s \in [0,1]$ such that $v(\s)=||v|| $ and considering $[t_0, t_1]$ as either of  $[0, \s]$ and $[\s, 1]$, we have
$$
v(t)  \geq  t||v|| \quad {\rm for} \quad t \in [0, \s],\\[.2cm]
$$
and
$$
v(t)  \geq (1-t)||v|| \quad {\rm for} \quad t \in [\s,1].
$$
Hence,
$$
v(t) \geq  \min\{t, 1-t\}||v||, \quad t \in [0, 1].
$$
\epf

\slm\mlabel{lm-compact}
Assume \rm{(H1)} holds. Then $\vect{T} _{\l}(K) \subset K$ and $\vect{T}_{\l}: K \to K$ is compact and continuous.
\elm
\pf
Lemma ~\ref{lm2} implies that $\vect{T} _{\l}(K) \subset K.$  It is a standard procedure to prove  that $\vect{T}_{\l}: K \to K$ is compact and continuous.

Let $$ \G = \frac{1}{4}\int^{\frac{3}{4}}_{\frac{1}{4}}\Big(\int^{s}_{\frac{1}{4}} N \t^{N-1}d\t\Big)^{\frac{1}{N}}ds> 0.$$
\slm\mlabel{f_estimate_>*}
Assume \rm{(H1)} holds. Let $ \vect{v}=(v_1,...,v_n) \in K $ and $\eta > 0$. If there exists a component $f^i$ of $\vect{f}$ such that
$$f^i(\vect{v}(t)) \geq \var(\eta \sum_{i=1}^n v_i(t)) \for  t \in [\frac{1}{4}, \frac{3}{4}] $$
then
$$
\norm{\vect{T}_{\l}\vect{v}} \geq \var^{-1}(\l) \G \eta \norm{\vect{v}}.
$$
\elm
\pf Note, from the definition of $\vect{T}_{\l}\vect{v}$,  that $T_{\l}^i \vect{v}(0)$ is the maximum value of $T_{\l}^i \vect{v}$ on [0,1].
It follows that
 \begin{equation*}
 \begin{split}
 \norm{\vect{T}_{\l}\vect{v}} & \geq \sup_{t\in[0,1]} \abs{T^i_{\l}\vect{v}(t)} \\[.2cm]
 & \geq   \int^{\frac{3}{4}}_{\frac{1}{4}}\var^{-1}\Big(\int^s_{\frac{1}{4}} \l N \t^{N-1}f^i(\vect{v}(\t))d\t\Big)ds \\[.2cm]
 & \geq   \int^{\frac{3}{4}}_{\frac{1}{4}}\var^{-1}\Big(\int^{s}_{\frac{1}{4}} \l N \t^{N-1}\var(\eta \sum_{i=1}^n v_i(\t))d\t\Big)ds \\
 & \geq   \int^{\frac{3}{4}}_{\frac{1}{4}}\var^{-1}\Big(\int^{s}_{\frac{1}{4}} \l N \t^{N-1}\var( \frac{\eta}{4}\norm{\vect{v}})d\t\Big)ds \\
 & = \var^{-1}(\l) \G \eta \norm{\vect{v}}
 \end{split}
 \end{equation*}
\epf
For each $i=1,...,n$, define a new function $\hat{f}^i(t): \mathbb{R}_+ \to \mathbb{R}_+$ by
$$\hat{f}^i(t) =\max \{f^i(\vect{v}):\vect{v} \in \mathbb{R}_+^n \; \rm{and}\;  \norm{\vect{v}} \leq t \}.$$
Note that $\hat{f}^i_{0}=\lim_{t \to 0} \frac{\hat{f}^i(t)}{\var(t)}$ and $\hat{f}^i_{\infty}=\lim_{t \to \infty} \frac{\hat{f}^i(t)}{\var(t)}.$
The following results hold, which were proved in \cite{Wang2003}. For completeness we here give a proof.
\slm\cite{Wang2003}\mlabel{lm6}
Assume \rm{ (H1)} holds. Then $\hat{f}^i_{0}=f^i_{0}\; \text{and}\; \hat{f}^i_{\infty}=f^i_{\infty},$ $i=1,...,n.$
\elm
\pf It is easy to see that $\hat{f}^i_{0}=f^i_{0}$. For the second part, we consider the two cases, (a) $f^i(\vect{v})$ is bounded and (b) $f^i(\vect{v})$ is unbounded. For case (a),
it follows, from $\lim_{t \to \infty}\var(t)=\infty$, that
$\hat{f}^i_{\infty}=0=f^i_{\infty}$.
\noindent For  case (b), for any $\delta > 0$, let $M^i=\hat{f}^i(\delta)$
and
$$ N_{\delta}^i = \inf\{\norm{\vect{v}}: \vect{v} \in \mathbb{R}_+^n, \;\norm{\vect{v}} \geq \delta, f^i(\vect{v}) \geq M^i \} \geq \delta,
$$
then
$$
\max\{f^i(\vect{v}): \norm{\vect{v}} \leq N_{\delta}^i, \; \vect{v} \in \mathbb{R}_+^n \}=M^i=\max\{f^i(\vect{v}): \norm{\vect{v}} = N_{\delta}^i, \; \vect{v} \in \mathbb{R}_+^n \}.
$$
Thus, for any $\delta > 0$, there exists a $ N_{\delta}^i \geq \delta$ such that
$$ \hat{f}^i(t) = \max\{f^i(\vect{v}): N_{\delta}^i \leq \norm{\vect{v}} \leq t , \; \vect{v} \in \mathbb{R}_+^n \} \;\; {\rm for} \;\; t > N_{\delta}^i.
$$
Now, suppose that $f^i_{\infty} < \infty$. In other words, for any $\e>0$, there is a $\delta>0$ such that
\begin{equation}\label{e-d-defintion}
f^i_{\infty} -\e <\frac{f^i(\vect{v})}{\var(\norm{\vect{v}})} < f^i_{\infty} +\e,  \textrm{  for }  \vect{v} \in \mathbb{R}_+^n,\; \norm{\vect{v}} > \delta.
\end{equation}
Thus, for $t>N_{\delta}^i$, there exist $\vect{v}_1, \vect{v}_2 \in \mathbb{R}_+^n$ such that  $\norm{\vect{v}_1}=t$, $t \geq \norm{\vect{v}_2} \geq N_{\delta}^i$  and $f^i(\vect{v}_2)=\hat{f}^i(t)$. Therefore,
\begin{equation}\label{e-d-inq}
\frac{f^i(\vect{v}_1)}{\var(\norm{\vect{v}_1})} \leq  \frac{\hat{f}^i(t)}{\var(t)} = \frac{f^i(\vect{v}_2)}{\var(t)} \leq \frac{f^i(\vect{v}_2)}{\var(\norm{\vect{v}_2})}.
\end{equation}
(\ref{e-d-defintion}) and (\ref{e-d-inq}) yield that
\begin{equation}\label{e-d-new}
f^i_{\infty} -\e <\frac{\hat{f}^i(t)}{\var(t)} < f^i_{\infty} +\e \textrm{ for }  t>N_{\delta}^i.
\end{equation}
Hence $\hat{f}^i_{\infty}=f^i_{\infty}$. Similarly, we can show $\hat{f}^i_{\infty}=f^i_{\infty}$ if $f^i_{\infty}=\infty$.
\slm\mlabel{f_estimate_<*}
Assume \rm{\rm{(H1)}} holds and let $ r >0 $. If there exits an $\e > 0$ such that
$$
\hat{f}^i(r) \leq \var(\e r),\;\; i=1,...,n,
$$
then
$$
\norm{\vect{T}_{\l}\vect{v}} \leq \var^{-1}(\l) \e n\norm{\vect{v}} \;\; {\rm for} \;\; \vect{v} \in \partial\O_{r}.
$$
\elm
\pf From the definition of $T_{\l}$, for $\vect{v} \in \partial\O_{r}$, we have
 \begin{eqnarray*}
 \norm{\vect{T}_{\l}\vect{v}} & =&  \sum_{i=1}^n \sup_{t\in[0,1]} \abs{T^i_{\l}\vect{v}(t)}\\
  & \leq &   \sum_{i=1}^n \var^{-1}(\int^1_0 \l N \t^{N-1} f^i(\vect{v}(\t))d\t) \\
  & \leq & \sum_{i=1}^n \var^{-1}(\int^1_0  N \t^{N-1}d\t \l\hat{f}^i(r))\\
  & \leq & \sum_{i=1}^n \var^{-1}(\int^1_0  N \t^{N-1}d\t \l \var(\e r))\\
  & = & \var^{-1}(\l) \e n \norm{\vect{v}}.
  \end{eqnarray*}
\epf
The following two lemmas are weak forms of Lemmas \ref{f_estimate_>*} and \ref{f_estimate_<*}.
\slm\mlabel{f_estimate_>*_weak}
Assume \rm{(H1)-(H2)} hold. If $ \vect{v} \in \partial \O_{r}$, $r >0$, then
$$
\norm{\vect{T}_{\l}\vect{v}}  \geq  4 \var^{-1}(\l) \G \var^{-1}(\hat{m}_r) $$
where $\hat{m}_r=\min \{ f^i(\vect{v}):  \vect{v} \in \mathbb{R}_+^n \; \rm{and}\; \frac{r}{4} \leq \norm{\vect{v}} \leq r,$ $i=1,...,n\}>0.$
\elm
\pf Since $ f_{i}(\vect{v}(t)) \geq  \hat{m}_r =\var(\var^{-1}(\hat{m}_r))\; \rm{for}\; t \in [\frac{1}{4}, \frac{3}{4}],$
$i=1,...,n,$ it is easy to see that this lemma can be shown in a similar manner as in Lemma \ref{f_estimate_>*}.
\epf
\slm\mlabel{f_estimate_<*_weak}
Assume \rm{\rm{(H1)-(H2)}} hold. If $ \vect{v} \in \partial \O_{r}$, $r >0$, then
$$
\norm{\vect{T}_{\l}\vect{v}} \leq \var^{-1}(\l) \var^{-1}(\hat{M}_r) n,
$$
where $\hat{M}_r=\max\{ f^i(\vect{v}):  \vect{v} \in \mathbb{R}_+^n \; \rm{and}\; \norm{\vect{v}} \leq r, \,\, i=1,...,n\}>0$ and $n$ is the positive constant defined in Lemma \ref{f_estimate_<*}
\elm
\pf Since $ f_{i}(\vect{v}(t)) \leq \hat{M}_r = \var(\var^{-1}( \hat{M}_r)) \; \rm{for}\; t \in [0,1]$,
$i=1,...,n,$ it is easy to see that this lemma can be shown in a similar manner as in Lemma \ref{f_estimate_<*}.
\epf

\section{Proof of Theorem \ref{th1}}
\pf
Part (a). $\vect{f}_0=0$ implies that $f^i_0=0$, $i=1,...,n$. It follows from Lemma ~\ref{lm6} that $\hat{f}^i_0=0$, $i=1,...,n.$ Therefore, we can choose $r_1 > 0$
so that $\hat{f}^i(r_1) \le \var(\e ) \var(r_1),$ $i=1,...,n$, where the constant $\e> 0$ satisfies
$$
\var^{-1}(\l) \e n < 1.
$$
We have by Lemma ~\ref{f_estimate_<*} that
$$
\norm{\vect{T}_{\l}\vect{v}} \leq \var^{-1}(\l) \e n\norm{\vect{v}} < \norm{\vect{v}} \quad \textrm{for} \quad  \vect{v} \in \partial\O_{r_1}.
$$
Now, since $\vect{f}_{\infty} = \infty$, there exists a component $f^i$ of $\vect{f}$ such that $f^i_{\infty}=\infty$. Therefore, there is
an $\hat{H} > 0$ such that
$$
f^i(\vect{v}) \geq \var(\eta)\var(\norm{\vect{v}})
$$
for $ \vect{v}=(v_1,...,v_n) \in \mathbb{R}_+^n$ and $\norm{\vect{v}} \geq \hat{H}$ ,
where $\eta > 0$ is chosen so that
$$
\var^{-1}(\l) \G \eta > 1.
$$
Let $r_2 = \max\{2r_1,4\hat{H} \}$. If $ \vect{v}=(v_1,...,v_n) \in \partial \O_{r_2}$, then
$$ \min_{\frac{1}{4}\leq t \leq \frac{3}{4}} \sum_{i=1}^n v_i(t) \geq \frac{1}{4}
\norm{\vect{v}}= \frac{1}{4} r_2 \geq \hat{H},$$
which implies that
$$f^i(\vect{v}(t)) \geq \var(\eta)\var(\sum_{i=1}^n v_i(t)) = \var(\eta\sum_{i=1}^n v_i(t))\; \rm{for} \; t \in [\frac{1}{4}, \frac{3}{4}].
$$
It follows from Lemma ~\ref{f_estimate_>*} that
$$
\norm{\vect{T}_{\l}\vect{v}} \geq \var^{-1}(\l) \G \eta \norm{\vect{v}} > \norm{\vect{v}} \quad \textrm{for}\quad  \vect{v} \in \partial\O_{r_2}.
$$
By Lemma ~\ref{lm1},
$\vect{T}_{\l}$ has a fixed point $\vect{v} \in  \O_{r_2} \setminus \bar{\O}_{r_1}$, which is  the desired positive solution of (\ref{eq4}).

Part (b). If $\vect{f}_0 = \infty$, there exists a component $f^i$ such that $f^i_{0}=\infty$. Therefore,
there is an $r_1 > 0$ such that
$$
f^i(\vect{v}) \geq \var(\eta)\var(\norm{\vect{v}})
$$
for $ \vect{v}=(v_1,...,v_n) \in \mathbb{R}_+^n$ and $ \norm{\vect{v}} \leq r_1,$
where $\eta > 0$ is chosen so that
$$
\var^{-1}(\l) \G \eta > 1.
$$
If $ \vect{v}=(v_1,...,v_n) \in  \partial \O_{r_1}$, then
$$f^i(\vect{v}(t)) \geq \var(\eta)\var(\sum_{i=1}^n v_i(t)) = \var(\eta\sum_{i=1}^n v_i(t)), \;\; {\rm for } \;\; t \in [0,1].$$
Lemma ~\ref{f_estimate_>*} implies that
$$
\norm{\vect{T}_{\l}\vect{v}} \geq \var^{-1}(\l) \G \eta \norm{\vect{v}} > \norm{\vect{v}} \quad \textrm{for}\quad  \vect{v} \in \partial\O_{r_1}.
$$
We now determine $\O_{r_2}$.  $\vect{f}_{\infty}=0$ implies that $f^i_{\infty}=0$, $i=1,...,n$. It follows from Lemma ~\ref{lm6} that $\hat{f}^i_{\infty}=0$, $i=1,...,n.$
Therefore there is an $r_2>2r_1$ such that
$$
\hat{f}^i(r_2) \le \var(\e) \var(r_2),\;i=1,...,n,
$$
where the constant $\e > 0$ satisfies
$$
\var^{-1}(\l) \e n < 1.
$$
Thus, we have by Lemma ~\ref{f_estimate_<*} that
$$
\norm{\vect{T}_{\l}\vect{v}} \leq \var^{-1}(\l) \e n \norm{\vect{v}} < \norm{\vect{v}} \quad \textrm{for}\quad  \vect{v} \in \partial\O_{r_2}.
$$
By Lemma ~\ref{lm1}, $\vect{T}_{\l}$ has a fixed point in  $\O_{r_2} \setminus \bar{\O}_{r_1}$, which is the desired positive solution of (\ref{eq1}).
\epf
\section{Proof of Theorem \ref{th2}}
\pf
Part (a).
Fix a number $r_1 > 0$. Lemma \ref{f_estimate_>*_weak} implies that there exists a $\l_0 >0$ such that
$$
\norm{\vect{T}_{\l}\vect{v}} > \norm{\vect{v}}=r_1, \;\; \rm{for}\;\;\vect{v} \in  \partial \O_{r_1}, \l > \l_0.
$$
If $\vect{f}_0=0$, then $f^i_0=0$, $i=1,...,n$. It follows from Lemma \ref{lm6} that
$$
\hat{f}^i_0=0,\;i=1,...,n.
$$
Therefore, we can choose $0 < r_2 < r_1 $ so that
$$
\hat{f}^i(r_2) \le \var(\e ) \var(r_2), \;i=1,...,n,
$$
where the constant $\e> 0$ satisfies
$$
\var^{-1}(\l) \e n < 1.
$$
We have by Lemma ~\ref{f_estimate_<*} that
$$
\norm{\vect{T}_{\l}\vect{v}} \leq \var^{-1}(\l) \e n \norm{\vect{v}} < \norm{\vect{v}} \quad \textrm{for} \quad  \vect{v} \in \partial\O_{r_2}.
$$
If $\vect{f}_{\infty}=0$, then $f^i_{\infty}=0$, $i=1,...,n$. It follows from Lemma \ref{lm6} that $\hat{f}^i_{\infty}=0$, $i=1,...,n.$
Therefore there is an $r_3>2r_1$ such that
$$
\hat{f}^i(r_3) \le \var(\e) \var(r_3),\;i=1,...,n,
$$
where the constant $\e > 0$ satisfies
$$
\var^{-1}(\l) \e n < 1.
$$
Thus, we have by Lemma ~\ref{f_estimate_<*} that
$$
\norm{\vect{T}_{\l}\vect{v}} \leq \var^{-1}(\l) \e n \norm{\vect{v}} < \norm{\vect{v}} \quad \textrm{for} \quad  \vect{v} \in \partial\O_{r_3}.
$$
It follows from Lemma ~\ref{lm1} that $\vect{T}_{\l}$ has a fixed point in  $\O_{r_1} \setminus \bar{\O}_{r_2}$ or $\O_{r_3} \setminus \bar{\O}_{r_1}$
according to $\vect{f}_0=0$
or $\vect{f}_{\infty}=0$, respectively. Consequently, (\ref{eq4}) has a positive solution for $ \l > \l_0$.

Part (b).
Fix a number $r_1 > 0$. Lemma \ref{f_estimate_<*_weak} implies that there exists a $\l_0 >0$ such that
$$
\norm{\vect{T}_{\l}\vect{v}}  <  \norm{\vect{v}}=r_1, \; {\rm for} \; \vect{v} \in  \partial \O_{r_1},\; 0< \l < \l_0.
$$
If $\vect{f}_0 = \infty$, there exists a component $f^i$ of $\vect{f}$ such that $f^i_{0}=\infty$. Therefore,
there is a positive number $r_2 <  r_1$ such that
$$
f^i(\vect{v}) \geq \var(\eta)\var(\norm{\vect{v}})
$$
for $ \vect{v}=(v_1,...,v_n) \in \mathbb{R}_+^n$ and $ \norm{\vect{v}} \leq r_2,$
where $\eta > 0$ is chosen so that
$$
\var^{-1}(\l) \G \eta > 1.
$$
Then
$$f^i(\vect{v}(t)) \geq \var(\eta)\var(\sum_{i=1}^n v_i(t)) = \var(\eta\sum_{i=1}^n v_i(t)),$$
for $\vect{v}=(v_1,...,v_n) \in  \partial \O_{r_2}, \;\; t \in [0,1].$
Lemma ~\ref{f_estimate_>*} implies that
$$
\norm{\vect{T}_{\l}\vect{v}} \geq \var^{-1}(\l) \G \eta \norm{\vect{v}} > \norm{\vect{v}} \quad \textrm{for} \quad  \vect{v} \in \partial\O_{r_2}.
$$
If $\vect{f}_{\infty} = \infty$, there exists a component $f^i$ of $\vect{f}$ such that $f^i_{\infty}=\infty$.  Therefore, there is
an $\hat{H} > 0$ such that
$$
f^i(\vect{v}) \geq \var(\eta)\var(\norm{\vect{v}})
$$
for $ \vect{v}=(v_1,...,v_n) \in \mathbb{R}_+^n$ and $\norm{\vect{v}} \geq \hat{H}$ ,
where $\eta > 0$ is chosen so that
$$
\var^{-1}(\l) \G \eta > 1.
$$
Let $r_3 = \max\{2r_1,4\hat{H} \}$. If $ \vect{v}=(v_1,...,v_n) \in \partial \O_{r_3}$, then
$$ \min_{\frac{1}{4}\leq t \leq \frac{3}{4}} \sum_{i=1}^n v_i(t) \geq \frac{1}{4}
\norm{\vect{v}}= \frac{1}{4} r_3 \geq \hat{H},$$
which implies that
$$
f^i(\vect{v}(t)) \geq \var(\eta)\var(\sum_{i=1}^n v_i(t)) = \var(\eta\sum_{i=1}^n v_i(t))\; \rm{for} \; t \in [\frac{1}{4}, \frac{3}{4}].
$$
It follows from Lemma ~\ref{f_estimate_>*} that
$$
\norm{\vect{T}_{\l}\vect{v}} \geq \var^{-1}(\l)\G \eta \norm{\vect{v}} > \norm{\vect{v}} \quad \rm{for} \quad  \vect{v} \in \partial\O_{r_3}.
$$
It follows from Lemma ~\ref{lm1} that
$\vect{T}_{\l}$ has a fixed point in  $\O_{r_1} \setminus \bar{\O}_{r_2}$ or $\O_{r_3} \setminus \bar{\O}_{r_1}$ according to $f_0=\infty$
or $f_{\infty}=\infty$, respectively. Consequently, (\ref{eq4}) has a positive solution for $ 0< \l < \l_0$.

Part (c).
Fix two numbers $0 < r_{3} < r_{4}.$ Lemma \ref{f_estimate_>*_weak} implies that there exists a $\l_0 >0$ such that for $\l >  \l_{0}$,
$$
\norm{\vect{T}_{\l}\vect{v}} > \norm{\vect{v}}, \;\;\rm{for} \;\; \vect{v} \in  \partial
\O_{r_i}, \;\; (i=3,4).
$$
Since $\vect{f}_0=0$ and ${\vect{f}_{\infty}=0},$ it follows from the proof of Theorem ~\ref{th2} (a) that we can choose $0< r_1 < r_3/2$ and $r_2>2r_4$  such that
$$\norm{\vect{T}_{\l}\vect{v}}  <  \norm{\vect{v}}, \;\; {\rm for} \;\; \vect{v} \in \partial \O_{r_i}, \;\; (i=1,2).$$
It follows from Lemma ~\ref{lm1} that
 $\vect{T}_{\l}$ has two fixed points $\vect{v}_1(t)$ and $\vect{v}_2(t)$ such that  $\vect{v}_1(t) \in \O_{r_3} \setminus \bar{\O}_{r_1}$ and $\vect{v}_2(t) \in \O_{r_2} \setminus \bar{\O}_{r_4}$ , which are
the desired distinct positive solutions of (\ref{eq4}) for $ \l > \l_0$ satisfying
$$
r_1 < \norm{\vect{v}_1} < r_3 < r_4 < \norm{\vect{v}_2} < r_2.
$$

Part (d).
Fix two numbers $0 < r_{3} < r_{4}.$ Lemma \ref{f_estimate_<*_weak} implies that there exists a $\l_0 >0$ such that for $ 0< \l < \l_{0} $,
$$
\norm{\vect{T}_{\l}\vect{v}} < \norm{\vect{v}}, \;\;\rm{for} \;\; \vect{v} \in  \partial
\O_{r_i}, \;\; (i=3,4).
$$
Since $\vect{f}_0=\infty$ and $\vect{f}_\infty=\infty$, it follows from the proof of Theorem ~\ref{th2} (b) that
we can choose $0< r_1 < r_3/2$ and $r_2>2r_4$ such that
$$\norm{\vect{T}_{\l}\vect{v}} > \norm{\vect{v}}, \;\; {\rm for} \;\; \vect{v} \in  \partial \O_{r_i}, \;\; (i=1,2).$$
It follows from Lemma ~\ref{lm1} that
$\vect{T}_{\l}$ has two fixed points $\vect{v}_1(t)$ and $\vect{v}_2(t)$ such that  $\vect{v}_1(t) \in \O_{r_3} \setminus \bar{\O}_{r_1}$ and $\vect{v}_2(t) \in \O_{r_2} \setminus \bar{\O}_{r_4}$ ,
which are the desired distinct positive solutions of (\ref{eq4}) for $ \l < \l_0$ satisfying
$$
r_1 < \norm{\vect{v}_1} < r_3 < r_4 < \norm{\vect{v}_2} < r_2.
$$

Part (e). Since $\vect{f}_0  < \infty$ and $\vect{f}_{\infty} < \infty$, then $f^i_{0}<\infty$ and $f^i_{\infty}<\infty$, $i=1,...,n.$
Therefore, for each $i=1,...,n,$ there exist positive numbers
$\e_1^i$, $ \e_2^i$, $r_1^i$ and $r_2^i$ such that $r_1^i < r_2^i$,
$$
  f^i(\vect{v})  \leq   \e_1^i \var(\norm{\vect{v}})\; \rm{for} \;\vect{v}\in \mathbb{R}_+^n,\; \norm{\vect{v}} \leq r_1^i,
$$
and
$$
  f^i(\vect{v})   \leq  \e_2^i \var(\norm{\vect{v}})\; \rm{for} \;\vect{v}\in \mathbb{R}_+^n,\; \norm{\vect{v}} \geq r_2^i.
$$
Let
$$
\e^i = \max\{\e_1^i, \e_2^i, \max\{\frac{f^i(\vect{v})}{\var(\norm{\vect{v}})}: \vect{v}\in \mathbb{R}_+^n, \;\; r_1^i \leq \norm{\vect{v}} \leq  r_2^i \}\} > 0
$$
and $ \e = \max\limits_{i=1,...,n} \{\e^i\} > 0.$ Thus, we have
$$
 f^i(\vect{v}) \leq  \e \var(\norm{\vect{v}})\; \rm{ for } \;\vect{v}\in \mathbb{R}_+^n,\;\;i=1,...,n.
$$
Assume $\vect{v}(t)$ is a positive solution of (\ref{eq4}). We will show that this leads to a contradiction
for $0<  \l < \l_0,$
where
$$
\l_0=\var(\frac{1}{n \var^{-1}(\e )}).
$$
In fact, for $0<  \l < \l_0$, since  $\vect{T}_{\l}\vect{v}(t) = \vect{v}(t)$ for $ t \in [0,1]$, we have
 \begin{eqnarray*}
  \norm{\vect{v}} & = &  \norm{\vect{T}_{\l}\vect{v}}\\
  & \leq &  \sum_{i=1}^n \var^{-1}(\int^1_0 N \t^{N-1}\e d\t \l\var(\norm{\vect{v}})) \\
  & = &  \sum_{i=1}^n \var^{-1}(\e \l\var(\norm{\vect{v}})) \\
  & = & \var^{-1}(\l)n\var^{-1}(\e) \norm{\vect{v}} \\
  & < & \norm{\vect{v}},
\end{eqnarray*}
which is a contradiction.

Part (f). Since $\vect{f}_0  > 0$ and $\vect{f}_{\infty} > 0$, there exist two components $f^i$ and $f^j$ of $\vect{f}$ such that
$f^i_0  > 0$ and $f^j_{\infty} > 0$. Therefore, there exist positive numbers
$\eta_1$, $ \eta_2$, $r_1$ and $r_2$ such that $r_1 < r_2,$
$$
  f^i(\vect{v})  \geq   \eta_1 \var(\norm{\vect{v}})\; \textrm{ for } \;\vect{v}\in \mathbb{R}_+^n,\; \norm{\vect{v}} \leq r_1,
$$
and
$$
f^j(\vect{v}) \geq  \eta_2 \var(\norm{\vect{v}})\; \textrm{ for } \;\vect{v}\in \mathbb{R}_+^n,\; \norm{\vect{v}} \geq r_2.
$$
Let
$$
\eta_3 = \min\{\eta_1, \eta_2, \min\{\frac{f^j(\vect{v})}{\var(\norm{\vect{v}})}: \vect{v}\in \mathbb{R}_+^n, \;\; \frac{r_1}{4} \leq \norm{\vect{v}} \leq  r_2\}\} > 0.
$$
Thus, we have
$$
  f^i(\vect{v}) \geq \eta_3 \var(\norm{\vect{v}})\; \textrm{ for } \;\vect{v}\in \mathbb{R}_+^n,\; \norm{\vect{v}} \leq r_1,
$$
and
$$
  f^j(\vect{v})   \geq   \eta_3 \var(\norm{\vect{v}})\; \textrm{ for } \;\vect{v}\in \mathbb{R}_+^n,\; \norm{\vect{v}} \geq \frac{r_1}{4}.
$$
Since $\eta_3\var(\norm{\vect{v}}) = \var\big(\var^{-1}(\eta_3)\big)\var(\norm{\vect{v}})$, it follows that

\begin{align} \label{ineq1}
  f^i(\vect{v})  & \geq   \var\big(\var^{-1}(\eta_3)\norm{\vect{v}}\big)\; \textrm{ for } \;\vect{v}\in \mathbb{R}_+^n,\; \norm{\vect{v}} \leq r_1
\end{align}
and
\begin{align} \label{ineq2}
  f^j(\vect{v})  & \geq  \var\big(\var^{-1}(\eta_3)\norm{\vect{v}}\big)\; \textrm{ for } \;\vect{v}\in \mathbb{R}_+^n,\; \norm{\vect{v}} \geq \frac{r_1}{4}.
\end{align}

Assume $\vect{v}(t)=(v_1,...,v_n)$ is a positive solution of (\ref{eq4}). We will show that this leads to a contradiction
for $\l > \l_0 =\var(\frac{1}{\G\var^{-1}(\eta_3)})$.  In fact, if $\norm{\vect{v}} \leq r_1$, (\ref{ineq1}) implies that
$$
f^i(\vect{v}(t)) \geq \var\big(\var^{-1}(\eta_3)\sum_{i=1}^n v_i(t)\big), \;\; {\rm for }\;\; t \in [0,1].
$$
On the other hand, if $\norm{\vect{v}} > r_1$, then $ \min_{\frac{1}{4}\leq t \leq \frac{3}{4}} \sum_{i=1}^n v_i(t) \geq \frac{1}{4}
\norm{\vect{v}} > \frac{1}{4} r_1,$
which, together with (\ref{ineq2}), implies that
$$
f^j(\vect{v}(t)) \geq \var\big(\var^{-1}(\eta_3)\sum_{i=1}^n v_i(t)\big), \;\; {\rm for }\;\; t \in [\frac{1}{4}, \frac{3}{4}].
$$
Since  $\vect{T}_{\l}\vect{v}(t) = \vect{v}(t)$ for $ t \in [0,1]$, it
follows from Lemma ~\ref{f_estimate_>*} that, for $\l > \l_0$,
 \begin{eqnarray*}
  \norm{\vect{v}} & = & \norm{\vect{T}_{\l}\vect{v}} \\[.2cm]
  & \geq &  \var^{-1}(\l) \G \var^{-1}(\eta_3) \norm{\vect{v}} \\[.2cm]
  & > & \norm{\vect{v}},
\end{eqnarray*}
which is a contradiction.\epf

\end{document}